\documentclass[a4paper,12pt,reqno]{amsart}
\usepackage{amsmath}
\usepackage{amssymb}
\usepackage{amsthm}

\usepackage{amscd}
\usepackage{amsfonts}
\usepackage{setspace}
\usepackage{version}

\title[Smooth numbers in arithmetic progressions]{On a paper of K. Soundararajan on smooth numbers in arithmetic progressions}
\author{Adam J Harper}
\address{Department of Pure Mathematics and Mathematical Statistics, Wilberforce Road, Cambridge CB\textup{3} \textup{0}WA, England}
\email{A.J.Harper@dpmms.cam.ac.uk}
\date{10th March 2011}
\thanks{The author is supported by a studentship from the Engineering and Physical Sciences Research Council of the United Kingdom.}

\numberwithin{equation}{section}
\theoremstyle{plain}

\onehalfspace
\newcommand{\Z}{\mathbb{Z}}
\newcommand{\C}{\mathbb{C}}

\newtheorem{conj1}{Conjecture}
\newtheorem{thm1}{Theorem}
\newtheorem{thm2}[thm1]{Theorem}
\newtheorem{prop1}{Proposition}
\newtheorem{lem1}{Lemma}
\newtheorem{majpr}{Majorant Principle}
\newtheorem{lem2}[lem1]{Lemma}
\newtheorem{rodbnd}{Rodosski\u{\i} Bound}
\newtheorem{rodbnd2}[rodbnd]{Rodosski\u{\i} Bound}
\newtheorem{charbnd}{Character Sum Bound}
\newtheorem{lem3}[lem1]{Lemma}

\begin{document}

\maketitle

\begin{abstract}
In a recent paper, K. Soundararajan showed, roughly speaking, that the integers smaller than $x$ whose prime factors are less than $y$ are asymptotically equidistributed in arithmetic progressions to modulus $q$, provided that $y^{4\sqrt{e} - \delta} \geq q$ and that $y$ is neither too large nor too small compared with $x$. We show that these latter restrictions on $y$ are unnecessary, thereby proving a conjecture of Soundararajan. Our argument uses a simple majorant principle for trigonometric sums to handle a saddle point that is close to 1.
\end{abstract}

\section{Introduction}
For $y \geq 1$, let $\mathcal{S}(y)$ denote the set of $y$-smooth numbers: that is, the set of numbers all of whose prime factors are less than or equal to $y$. For $x \geq 1$, and natural numbers $a,q$, we define the following counting functions:
$$ \Psi_{q}(x,y) := \sum_{n \leq x, (n,q)=1} \textbf{1}_{\{n \in \mathcal{S}(y)\}}, \;\;\; \textrm{ and } \;\;\; \Psi(x,y;q,a) := \sum_{n \leq x, n \equiv a (\textrm{mod } q)} \textbf{1}_{\{n \in \mathcal{S}(y)\}}, $$
where $\textbf{1}$ denotes the indicator function.

In his 2008 article~\cite{sound}, K. Soundararajan makes the following equidistribution conjecture:
\begin{conj1}[Soundararajan, 2008]
Let $A$ be a given positive real number. Let $y$ and $q$ be large with $q \leq y^{A}$, and let $(a,q)=1$. Then as $\log x/\log q \rightarrow \infty$ we have
$$ \Psi(x,y;q,a) \sim \frac{1}{\phi(q)}\Psi_{q}(x,y). $$
\end{conj1}
As Soundararajan~\cite{sound} discusses, in our current state of knowledge about character sums it would be very hard to prove the conjecture for $A \geq 4\sqrt{e}$; for if the conjecture is true, then e.g. the least quadratic non-residue modulo $q$ must lie below $q^{1/A}$. However, Soundararajan is able to prove the conjecture for $A < 4\sqrt{e}$, on the additional assumption that $e^{y^{1-\epsilon}} \geq x \geq y^{(\log\log y)^{4}}.$ In this note we establish the following result, confirming that this assumption on $y$ is not needed.
\begin{thm1}
Let $\delta > 0$, and suppose that $y \leq x$, that $2 \leq q \leq y^{4\sqrt{e} - \delta}$, and $(a,q)=1$. If $y$ is large enough depending on $\delta$, then as $\log x / \log q \rightarrow \infty$ we have
$$ \Psi(x,y;q,a) \sim \frac{1}{\phi(q)}\Psi_{q}(x,y). $$
\end{thm1}
In fact, following Soundararajan~\cite{sound}, we work with a smooth weight function $\Phi(n/x)$ throughout: see $\S 2.1$ for further details. We obtain a smoothly weighted version of Theorem 1 with a quantitative error term
$$ O_{\Phi}\left(\frac{\Psi_{q}(x,y)}{\phi(q)} \left(\min\{\frac{q\sqrt{y}}{(\log\log x)^{1/3}}, \frac{1}{\log y}\} + \frac{\log q}{u^{c} \log y} + \frac{\log w}{(c \delta w^{\delta/2})^{c \log(2+(\log x)/ y)}} \right) \right) , $$
where $c > 0$ is an absolute constant, and we write $u=\log x / \log y$, $v=\log x /\log q$ and $w = \min\{v,y\}$. We caution that the reader should not try to read off the necessary dependence of $y$ on $\delta$ from this bound, as it is not valid unless $y$ is large enough that the error term $O_{\delta}(y^{-\delta^{2}}\log^{2}y)$ in Character Sum Bound 1 (see $\S 2.5$, below) is smaller than $\delta/2$, say.

Soundararajan's article~\cite{sound} also contains an `equidistribution in cosets' result, which gives information towards Conjecture 1 in the case $A \geq 4\sqrt{e}$. This is again proved on the assumption that $e^{y^{1-\epsilon}} \geq x \geq y^{(\log\log y)^{4}}$, which we can now remove.
\begin{thm2}
Let $A$ be a given positive real number, and let $y$ and $q$ be large with $q \leq y^{A}$. There is a subgroup $H$ of $(\Z/q\Z)^{*}$, of index bounded in terms of $A$ only, such that whenever $a,b \in (\Z/q\Z)^{*}$ satisfy $a/b \in H$ we have
$$ \Psi(x,y;q,a) - \Psi(x,y;q,b) = o(\Psi_{q}(x,y)/\phi(q)) \;\;\;\;\; \textrm{ as } \log x / \log q \rightarrow \infty. $$
\end{thm2}
We will not say much about Theorem 2, which follows from the proof of Theorem 1 as in Soundararajan's paper~\cite{sound}. The author has tried to keep explicit dependence on $A$ in the arguments below, so a keen reader may check that provided $\log y / ((A+1)\log(A+2))$ and $\log v / ((A+1)\log(A+2))$ are sufficiently large we have a smoothly weighted version with a better error term than above (the term involving $w$ can be removed). The error term is better in Theorem 2 because we do not need Character Sum Bound 1 for the proof. In our first appendix we comment briefly on how to pass from the smoothly weighted to the unweighted version of Theorem 2, which seems to require a slightly different procedure\footnote{A reader who is checking Theorem 2 may wish to consult this appendix first. The unsmoothing procedure that we use will allow one to prove Theorem 2 without needing to analyse the characters $\chi \in \mathcal{B}$, (defined in $\S 2.5$, below), which is quite a helpful simplification.} now that $y$ is unrestricted in terms of $x$.

\vspace{12pt}
Soundararajan~\cite{sound} argues, roughly, by observing the usual decomposition
$$ \Psi(x,y;q,a) = \frac{1}{\phi(q)} \sum_{\chi (\textrm{mod } q)} \overline{\chi(a)} \Psi(x,y;\chi), \;\;\; \textrm{ where } \;\;\; \Psi(x,y;\chi) := \sum_{n \leq x} \chi(n)\textbf{1}_{\{n \in \mathcal{S}(y)\}},$$
and analysing $\Psi(x,y;\chi)$ using knowledge of the $L$-series $L(s,\chi)$. His key innovation, perhaps, is to exploit the fact that we are interested in all characters to modulus $q$ taken together, and that we can make much stronger statistical statements than we can statements about individual $L$-series.

Perhaps surprisingly, it is when $y$ is close to $x$ that Soundararajan's argument is difficult to extend. This is because of a `saddle point problem': as $y$ approaches $x$, the (Euler product) terms from which one can gain by making non-trivial estimations carry progressively less weight, so it is important not to lose anywhere else. To achieve this we avoid applying absolute value bounds to integrals, and instead exploit a majorant principle for trigonometric sums. See $\S\S 2.3-2.5$. It is the author's opinion that this argument is the most interesting new aspect of this work.

A further difficulty in establishing Theorem 1 is that two parts of Soundararajan's proof, his ``basic argument'' and ``Rodosski\u{\i} argument'', are valid respectively when a quantity $k$ (explained in $\S 2.1$) is quite large depending on $x$, or is quite small depending on $u$. When $y$ approaches $x$ a gap emerges between these ranges, and to deal with this we need an argument based on Taylor expansion and a smoothed explicit formula. See $\S 2.3$ and $\S 3$. When $y$ is small Soundararajan's proofs~\cite{sound} almost go through, except for minor technical problems and some difficulties if $y$ does not tend to infinity with $x$. In $\S\S 2.6-2.7$ we give an argument that addresses these problems.

\vspace{12pt}
The author has tried to write this note in a reasonably self-contained way, whilst not simply repeating arguments that appear in Soundararajan's paper~\cite{sound}. To this end, three important pieces of `$L$-function information' obtained by Soundararajan are stated without proof in $\S 2$, as Rodosski\u{\i} Bound 1, Rodosski\u{\i} Bound 2 and Character Sum Bound 1. Except in the application of these bounds, (and the general set-up, which we recall in $\S 2.1$), many details of our argument are different from that of Soundararajan~\cite{sound}, and so we give a detailed account.

Since it adds very little complication, and may be illuminating, we shall prove Theorem 1 for all $y$, and not only the range not covered by Soundararajan's results. We distinguish in our work between ``large $y$'', namely $e^{\log^{1/10}x} < y \leq x$; ``small $y$'', namely $(\log\log x)^{3} \leq y \leq e^{\log^{1/10}x}$; and ``very small $y$'', namely $y < (\log\log x)^{3}$.

If $q < \sqrt{y}$, say, a result of Granville~\cite{granville} shows $\Psi(x,y;q,a)$ is $\Psi_{q}(x,y)\phi(q)^{-1}(1 + O(\log^{-1}y (1+u^{-c}\log q)))$. We invoke this result except when $q\sqrt{y} \leq (\log\log x)^{1/3}$ (for which see $\S 2.7$), and can therefore always assume that $q \geq \sqrt{y}$ except in that case. This will be convenient in applying various $L$-function computations, so that $\log q$ is somewhat comparable to $\log y$. The reader should also bear in mind, when checking that our proof supplies the bound claimed, that if $\sqrt{y} \leq q \leq y^{4\sqrt{e}}$ then $u=\log x /\log y$ is comparable in size to $v = \log x /\log q$.

\section{Overview of the argument}

\subsection{Initial set-up}
This subsection records some preliminary observations, mostly from $\S 2$ of Soundararajan's paper~\cite{sound} (which may be consulted for a more detailed description).

Let $\Phi : [0,\infty) \rightarrow [0,1]$ be a function supported on $[0,2]$, which equals 1 on $[0,1/2]$, and which is nine times continuously differentiable (say). We set
$$ \Psi(x,y;q,a,\Phi) := \sum_{n \in \mathcal{S}(y), n \equiv a (\textrm{mod } q)} \Phi(n/x), $$
which has a decomposition into weighted character sums $\Psi(x,y;\chi,\Phi)$ as in the introduction. We will work to show that $\Psi(x,y;q,a,\Phi)$ is approximately equal for all $(a,q)=1$, and by choosing $\Phi$ to bound $\textbf{1}_{[0,1]}$ from above and then below (in a way explained further in our first appendix) this will imply Theorem 1.

We define a truncated Euler product corresponding to $y$-smooth numbers, viz.
$$ L(s,\chi;y):=\prod_{p \leq y} \left(1-\frac{\chi(p)}{p^{s}}\right)^{-1} = \sum_{n \in \mathcal{S}(y)} \frac{\chi(n)}{n^{s}}, \;\;\; \Re(s) > 0. $$
Then, as usual, we can represent $\Psi(x,y;\chi,\Phi)$ as a contour integral involving $L(s,\chi;y)$:
$$ \Psi(x,y;\chi,\Phi) = \frac{1}{2\pi i} \int_{c - i\infty}^{c + i\infty} L(s,\chi;y) x^{s} \breve{\Phi}(s) ds, \;\;\; c > 0, $$
where $\breve{\Phi}(s) = \int_{0}^{\infty} \Phi(t) t^{s-1} dt$ is the Mellin transform of $\Phi$. Because $\Phi$ is so smooth, integration by parts shows $|\breve{\Phi}(s)| \ll_{\Phi} |s|^{-1}(|s|+1)^{-8}$ for $\Re(s) > 0$, which will be used several times. We also note that $\breve{\Phi}(c) \geq \int_{0}^{1/2} t^{c-1} dt \geq 1/(2c)$ for $0 < c \leq 1$.

We choose $c$ to be $\alpha=\alpha(x,y)$, a quantity coming from a saddle-point argument of Hildebrand and Tenenbaum~\cite{ht} (that was extended\footnote{Hildebrand and Tenenbaum~\cite{ht} studied $\Psi(x,y):=\sum_{n \leq x} \textbf{1}_{\{n \in \mathcal{S}(y)\}}$, and later de la Bret\`{e}che and Tenenbaum~\cite{dlbten} showed that some `obvious' adaptations of their results also hold for $\Psi_{q}(x,y)$, on a wide range of $q$. The bound that we record for $\Psi(x,y;\chi_{0},\Phi)$, involving the smoothing $\textbf{1}_{[0,1/2]} \leq \Phi \leq \textbf{1}_{[0,2]}$, is an easy consequence of Th\'{e}or\`{e}me 2.1 of de la Bret\`{e}che and Tenenbaum~\cite{dlbten} (e.g. because, as in our first appendix, $\Psi_{q}(x/2,y)$ is comparable in size to $\Psi_{q}(2x,y)$).} to treat $\Psi_{q}(x,y)$ by de la Bret\`{e}che and Tenenbaum~\cite{dlbten}). In practice this means the following: provided that $y$, $u=\log x/\log y$ and $\log^{2}y / \log q$ (say) are larger than certain absolute constants, as we assume throughout, we have
$$ \alpha(x,y) = \left\{ \begin{array}{ll}
     1 - \frac{\log(u\log u)}{\log y} + O(\frac{1}{\log y}) & \textrm{if } y > \log x   \\
     \Theta(\frac{y}{u \log^{2}y}) & \textrm{otherwise}
\end{array} \right. ,$$
and, writing $\chi_{0}$ for the principal Dirichlet character to modulus $q$, we have
$$ \Psi(x,y;\chi_{0},\Phi) \gg \frac{x^{\alpha} L(\alpha,\chi_{0};y) \breve{\Phi}(\alpha)}{\sqrt{2\pi (1+\log x/y)\log x \log y}}. $$

Finally we present some notation concerning the zeros of the $L$-series $L(s,\chi)$. For $0 \leq k \leq (\log q)/2$, write
$$ \Xi(k) := \{\chi : \chi \neq \chi_{0}, L(\sigma+it,\chi) \neq 0 \textrm{ for } \sigma > 1-\frac{k}{\log q}, |t| \leq q, \textrm{ but } L(\sigma+it,\chi)=0 $$
$$ \textrm{ for some } \sigma > 1-\frac{k+1}{\log q}, |t| \leq q\}.$$
As Soundararajan~\cite{sound} describes, the so-called log-free zero density estimate implies that $\#\Xi(k) \leq C_{1}e^{C_{2}k}$ for all $k$, for certain absolute constants $C_{1},C_{2}$. Thus
$$ \Psi(x,y;q,a,\Phi) = \frac{\Psi(x,y;\chi_{0},\Phi)}{\phi(q)} + O(\frac{1}{\phi(q)} \sum_{0 \leq k \leq \log q/2} e^{C_{2}k} \max_{\chi \in \Xi(k)} \left|\int_{\alpha - i\infty}^{\alpha + i\infty} L(s,\chi;y) x^{s} \breve{\Phi}(s) ds \right|). $$

To prove Theorem 1, we will show that the ``big Oh'' term in the preceding equation is of smaller order than our lower bound for the main term. To this end we will separate the summation over $k$ into summations over different ranges, as described in $\S 2.2$. We see immediately, however, that
$$ \left|\int_{\alpha + i(yq)^{1/4}}^{\alpha + i\infty} L(s,\chi;y) x^{s} \breve{\Phi}(s) ds \right| \leq L(\alpha,\chi_{0};y) x^{\alpha} \int_{\alpha + i(yq)^{1/4}}^{\alpha + i\infty} |\breve{\Phi}(s)| ds \ll \frac{L(\alpha,\chi_{0};y) x^{\alpha} \breve{\Phi}(\alpha)}{\breve{\Phi}(\alpha) y^{2}q^{2}}, $$
because of the rapid decay of $\breve{\Phi}(s)$. This is of smaller order than $\Psi(x,y;\chi_{0},\Phi)/yq^{2}$, (using the lower bound $\breve{\Phi}(\alpha) \geq 1/(2\alpha) \gg (\log x \log y )/y$ if $y \leq \log x$), and clearly the same holds for the integral over $(\alpha-i\infty,\alpha-i(yq)^{1/4}]$. Thus, unless $y$ is very small (for which see $\S 2.7$), it will suffice to prove satisfactory bounds for $\int_{\alpha - i(yq)^{1/4}}^{\alpha + i(yq)^{1/4}} L(s,\chi;y) x^{s} \breve{\Phi}(s) ds$, for $\chi \in \Xi(k)$. Note that $(yq)^{1/4} \leq q^{3/4}$ on our assumption that $q \geq \sqrt{y}$.

\subsection{Ranges of $k$}
For ``large'' $y$, in the sense of the introduction, we separate the summation over $0 \leq k \leq (\log q)/2$ into three different ranges, as follows:
\begin{itemize}
\item the ``basic range'', $\sqrt{u} \leq k \leq (\log q)/2$;

\item the ``Rodosski\u{\i} range'', $4A\log A + D \leq k < \sqrt{u}$, where $D$ is the absolute constant appearing in Rodosski\u{\i} Bound 1 in $\S 2.4$;

\item the ``problem range'', $0 \leq k < 4A\log A + D$.
\end{itemize}
(In Theorem 1 we have $A = 4\sqrt{e}-\delta$, and the reader may think of $A$ simply as $O(1)$. We continue to explicitly record dependence on $A$ to aid anyone checking Theorem 2.)

This is analogous to Soundararajan's argument~\cite{sound}, but our definitions of the ranges are different. In $\S\S 2.3-2.5$ we study these ranges in turn, and the reader may compare with $\S\S 3-5$ of Soundararajan's article~\cite{sound}.

For smaller $y$ the situation is simpler because one can treat the ``basic range'' and the ``Rodosski\u{\i} range'' in a unified way. This is discussed in $\S\S 2.6-2.7$.

\subsection{A modified zero-free region argument}
In $\S 3$ we will prove the following result, which we will need in place of Lemma 3.2 of Soundararajan~\cite{sound}:
\begin{prop1}
Let $B > 0$ be fixed, and let $y \geq 2$ and $\sqrt{y} \leq q \leq y^{A}$. If $\log x / ((B+1)\log q)$ and $k/((A+1)\log(u\log u))$ are larger than certain absolute constants, and if $y$ is ``large'' then the following holds. For any $\chi \in \Xi(k)$, any $\alpha-Bk/\log x \leq \sigma \leq \alpha$, and any $|t| \leq q/2$, we have
$$ |\log L(\sigma+it,\chi;y) - \log L(\alpha+it,\chi;y)| \leq k/50+O(1). $$
\end{prop1}
The proof of this involves using a smoothed explicit formula to analyse the first and second derivatives of $\log L(\sigma+it,\chi;y)$.

We will also use the following consequence of Fubini's theorem and the Cauchy--Schwarz inequality, whose proof is an easy exercise:
\begin{lem1}
Suppose that $\beta, r > 0$, and that $F(s)$ is any integrable function on the interval $[\beta,\beta+ir] \subseteq \C$. Let $G(s)$ be an Euler product of the following form:
$$G(s) := \prod_{p \leq y} \left(1-\frac{g(p)}{p^{s}}\right)^{-1}, $$
where $y \geq 0$ is fixed, and $|g(p)| \leq 1$ for all primes $p$.
Then
$$ \left|\int_{\beta}^{\beta+ir} G(s) F(s) ds \right| \leq M (|G(\beta)| + \sqrt{\int_{\beta}^{\beta+ir} |\frac{G'}{G}(s)|^{2} d|s| \int_{\beta}^{\beta+ir} |G(s)|^{2} d|s|} ), $$
where
$$ M = M(\beta,r,F) := \sup_{0 \leq t \leq r} \left|\int_{\beta+it}^{\beta+ir} F(s)ds \right| . $$
\end{lem1}

We apply Proposition 1 with $B$ chosen as $C_{2}+2$, where $C_{2}$ is the constant in the log-free zero density estimate in $\S 2.1$. Using this together with the rapid decay of $\breve{\Phi}(s)$, we note firstly that, {\em under the conditions of Proposition 1},
\begin{eqnarray}
&& \int_{\alpha - i(yq)^{1/4}}^{\alpha + i(yq)^{1/4}} L(s,\chi;y)x^{s}\breve{\Phi}(s) ds \nonumber \\
& = & \int_{\alpha - \frac{Bk}{\log x} -i\log^{1/4}y}^{\alpha - \frac{Bk}{\log x} +i\log^{1/4}y} L(s,\chi;y)x^{s}\breve{\Phi}(s) ds + \int_{\alpha - \frac{Bk}{\log x} -i(yq)^{1/4}}^{\alpha - \frac{Bk}{\log x} -i\log^{1/4}y} L(s,\chi;y)x^{s}\breve{\Phi}(s) ds + \nonumber \\
&& + \int_{\alpha - \frac{Bk}{\log x} +i\log^{1/4}y}^{\alpha - \frac{Bk}{\log x} +i(yq)^{1/4}} L(s,\chi;y)x^{s}\breve{\Phi}(s) ds + O(e^{k/50}L(\alpha,\chi_{0};y) x^{\alpha}/y^{2}q^{2}). \nonumber
\end{eqnarray}
The second and third integrals may also be estimated just using Proposition 1, showing they are $ O(e^{k/50} L(\alpha,\chi_{0};y) x^{\alpha-Bk/\log x}/\log^{2}y)$. Both ``big Oh'' terms are
$$ \ll L(\alpha,\chi_{0};y) x^{\alpha} (\frac{1}{y^{2} q^{199/100}} + \frac{e^{-(C_{2}+99/50)k} \sqrt{u}}{\sqrt{\log x} \log^{3/2}y} ) \ll \Psi(x,y; \chi_{0},\Phi)(\frac{1}{yq^{199/100}} + \frac{e^{-(C_{2}+1)k}}{\log y}), $$
on recalling our lower bound for $\Psi(x,y; \chi_{0},\Phi)$ and that $\log(u\log u) \ll k \leq (\log q)/2$.

Combining Lemma 1 (with the choices $F(s)=x^{s}\breve{\Phi}(s)$ and $G(s)=L(s,\chi;y)$) with Proposition 1, we see the first integral is
\begin{eqnarray}
& \ll & \sup_{0 \leq t \leq \log^{1/4}y} \left|\int_{\alpha - \frac{Bk}{\log x} +it}^{\alpha - \frac{Bk}{\log x} +i\log^{1/4}y} x^{s}\breve{\Phi}(s) ds \right| \cdot \biggl( |L(\alpha - \frac{Bk}{\log x},\chi;y)| + \nonumber \\
&& + \sqrt{\int_{\alpha - \frac{Bk}{\log x} - i\log^{1/4}y}^{\alpha - \frac{Bk}{\log x} + i\log^{1/4}y} |\frac{L'(s,\chi;y)}{L(s,\chi;y)}|^{2} d|s| \int_{\alpha - \frac{Bk}{\log x} - i\log^{1/4}y}^{\alpha - \frac{Bk}{\log x} + i\log^{1/4}y} |L(s,\chi;y)|^{2} d|s|} \biggr) \nonumber \\
& \ll & e^{-Bk+k/50} x^{\alpha} \sup_{0 \leq t \leq \log^{1/4}y} \left|\int_{t}^{\log^{1/4}y} x^{ir}\breve{\Phi}(\alpha-\frac{Bk}{\log x}+ir) dr \right| \nonumber \\
&& \cdot (L(\alpha,\chi_{0};y) + \sqrt{\int_{\alpha - \frac{Bk}{\log x} - i\log^{1/4}y}^{\alpha - \frac{Bk}{\log x} + i\log^{1/4}y} |\frac{L'(s,\chi;y)}{L(s,\chi;y)}|^{2} d|s| \int_{\alpha -i\log^{1/4}y}^{\alpha + i\log^{1/4}y} |L(s,\chi;y)|^{2} d|s|} ) . \nonumber
\end{eqnarray}

At this point we invoke the following majorant principle for trigonometric sums, which we quote from chapter 7.3 of Montgomery's book~\cite{mont}:
\begin{majpr}[Wirsing, and others]
Let $\lambda_{1},...,\lambda_{N}$ be real numbers, and suppose that $|a_{n}| \leq A_{n}$ for all $n$. Then
$$ \int_{-T}^{T} \left| \sum_{n=1}^{N} a_{n}e^{2\pi i\lambda_{n} t} \right|^{2} dt \leq 3 \int_{-T}^{T} \left| \sum_{n=1}^{N} A_{n}e^{2\pi i\lambda_{n} t} \right|^{2} dt. $$
\end{majpr}
The point is that $L(s,\chi;y)$ and its logarithmic derivative $L'(s,\chi;y)/L(s,\chi;y)$ are Dirichlet series, so in particular are trigonometric series with $\lambda_{n}$ chosen as $-\log n /2\pi$. (The majorant principle is stated for finite sums, but it remains valid for uniformly convergent Dirichlet series, as Montgomery~\cite{mont} remarks\footnote{See Chapter III.4.3 of Tenenbaum~\cite{ten} for an application of Majorant Principle 1 to Dirichlet series, concerning means of multiplicative functions. The author thanks K. Soundararajan for this reference.}.) Thus it will suffice to estimate the above with the squareroot term replaced by
$$ \sqrt{\int_{\alpha - \frac{Bk}{\log x} - i\log^{1/4}y}^{\alpha - \frac{Bk}{\log x} + i\log^{1/4}y} |\frac{L'(s,1;y)}{L(s,1;y)}|^{2} d|s| \int_{\alpha -i\log^{1/4}y}^{\alpha + i\log^{1/4}y} |L(s,\chi_{0};y)|^{2} d|s|}. $$

In our second appendix we show how to estimate the remaining integrals, which is fairly standard. It turns out, provided $y$ is ``large'' and $Bk/\log x \leq 1/8$, say, (so $\alpha - Bk/\log x \geq 3/4$), that the whole of the above is
$$ \ll e^{-Bk+k/50} x^{\alpha} \cdot \frac{L(\alpha,\chi_{0};y)}{\log x} (1 + \sqrt{y^{2Bk/\log x} u^{2} \log u}) \ll \Psi(x,y;\chi_{0},\Phi) \sqrt{u\log u} e^{-(C_{2}+99/50-B/u)k}. $$
Thus if $u/(A+1)^{3}$ is larger than an absolute constant, (so $\sqrt{u}/((A+1)\log(u\log u))$ is large, and therefore Proposition 1 is applicable for $k$ in the ``basic range''), we have the more than satisfactory estimate
$$\left|\int_{\alpha - i(yq)^{1/4}}^{\alpha + i(yq)^{1/4}} L(s,\chi;y) x^{s} \breve{\Phi}(s) ds \right| \ll \Psi(x,y; \chi_{0},\Phi)(\frac{1}{yq^{199/100}} + e^{-(C_{2}+1)k}) $$
when $k$ is in that range\footnote{We could use our argument on a much wider range of $k$ than the ``basic range''. However, this would not quite be large enough to dispense with the ``Rodosski\u{\i} range'' argument in $\S 2.4$, and splitting into ranges as we do yields better quantitative bounds on our integrals.}.

\subsection{A modified Rodosski\u{\i} argument}
We modify the ``Rodosski\u{\i} type argument'' from Soundararajan's paper~\cite{sound} (in which zeros of $L$-series are studied with carefully chosen weights) in the same kind of way as the zero-free region argument, by using Fubini's theorem and Majorant Principle 1. We begin with a variant of Lemma 1.
\begin{lem2}
Suppose that $\beta, r, F(s), G(s)$ are as in the statement of Lemma 1. Then
$$ \left|\int_{\beta}^{\beta+ir} G(s) F(s) ds \right| \leq M^{*} (|G^{*}(\beta)| + \sqrt{\int_{\beta}^{\beta+ir} |\frac{G'}{G}(s)|^{2} d|s| \int_{\beta}^{\beta+ir} |G^{*}(s)|^{2} d|s|} ), $$
where
$$G^{*}(s) := \prod_{p \leq \sqrt{y}} \left(1-\frac{g(p)}{p^{s}}\right)^{-1}, \; \textrm{ and } \; M^{*} := \sup_{0 \leq t \leq r} \left( \left|\int_{\beta+it}^{\beta+ir} F(s)ds \right| \prod_{\sqrt{y} < p \leq y} \left|1-\frac{g(p)}{p^{\beta+it}}\right|^{-1} \right) . $$
\end{lem2}

Similarly to $\S 2.3$, the rapid decay of $\breve{\Phi}(s)$ implies that
\begin{eqnarray}
&& \int_{\alpha - i(yq)^{1/4}}^{\alpha + i(yq)^{1/4}} L(s,\chi;y) x^{s} \breve{\Phi}(s) ds \nonumber \\
& = & \int_{\alpha - i\log^{1/4}y}^{\alpha + i\log^{1/4}y} L(s,\chi;y) x^{s} \breve{\Phi}(s) ds + O(\frac{x^{\alpha} L(\alpha,\chi_{0};\sqrt{y})}{\log^{2}y} \sup_{\log^{1/4}y < |t| \leq (yq)^{1/4}} \prod_{\sqrt{y} < p \leq y} \left|1-\frac{\chi(p)}{p^{\alpha+it}}\right|^{-1}). \nonumber
\end{eqnarray}
Combining Lemma 2 with Majorant Principle 1, for ``large'' $y$ we see
\begin{eqnarray}
&& \int_{\alpha - i\log^{1/4}y}^{\alpha + i\log^{1/4}y} L(s,\chi;y)x^{s}\breve{\Phi}(s) ds \nonumber \\
& \ll & \sup_{0 \leq t \leq \log^{1/4}y} \left|\int_{\alpha +it}^{\alpha +i\log^{1/4}y} x^{s}\breve{\Phi}(s) ds \right| \cdot \sup_{|t| \leq \log^{1/4}y} \prod_{\sqrt{y} < p \leq y} \left|1-\frac{\chi(p)}{p^{\alpha+it}}\right|^{-1} \nonumber \\
&& \cdot (L(\alpha,\chi_{0};\sqrt{y}) + \sqrt{\int_{\alpha - i\log^{1/4}y}^{\alpha + i\log^{1/4}y} |\frac{L'(s,1;y)}{L(s,1;y)}|^{2} d|s| \int_{\alpha - i\log^{1/4}y}^{\alpha + i\log^{1/4}y} |L(s,\chi_{0};\sqrt{y})|^{2} d|s|} ) \nonumber \\
& \ll & \frac{x^{\alpha}L(\alpha,\chi_{0};\sqrt{y}) u \sqrt{\log u}}{\log x} \sup_{|t| \leq \log^{1/4}y} \prod_{\sqrt{y} < p \leq y} \left|1-\frac{\chi(p)}{p^{\alpha+it}}\right|^{-1}. \nonumber
\end{eqnarray}
Here the final inequality again used the estimates from our second appendix.

Obtaining a non-trivial bound has now reduced to obtaining a sufficiently non-trivial estimate for the products over primes. This will follow from the next result, which is the content of Lemmas 4.2 and 4.3 of Soundararajan~\cite{sound}.
\begin{rodbnd}[Soundararajan, 2008]
There is an absolute constant $D$ for which the following is true. Suppose that $\chi \in \Xi(k)$ for some $k \geq 4A\log A + D$. If $q \leq y^{A}$, and $|t| \leq q/2$, and $y/(A+1)^{2}$ is large enough, then
$$ \sum_{\sqrt{y} \leq p \leq y, p \nmid q} \frac{1 - \Re(\chi(p)p^{-it})}{p} \log p \geq \frac{\log y}{5}. $$
\end{rodbnd}

Then
$$ \left|\frac{1-\chi(p)p^{-\alpha-it}}{1-p^{-\alpha}}\right| = \left|1 + \frac{1-\chi(p)p^{-it}}{p^{\alpha}-1}\right| \geq 1+\sum_{k=1}^{\infty} \frac{1-\Re(\chi(p)p^{-it})}{p^{k\alpha}} \geq e^{(1-\Re(\chi(p)p^{-it}))/p^{\alpha}}, $$
so (as in the argument of Lemma 4.2 of Soundararajan~\cite{sound}) we find
\begin{eqnarray}
\frac{L(\alpha,\chi_{0};\sqrt{y})}{L(\alpha,\chi_{0};y)} \sup_{|t| \leq (yq)^{1/4}} \prod_{\sqrt{y} < p \leq y} \left|1-\frac{\chi(p)}{p^{\alpha+it}}\right|^{-1} & = &  \sup_{|t| \leq (yq)^{1/4}} \prod_{\sqrt{y} < p \leq y} \left|\frac{1-\chi(p)p^{-\alpha-it}}{1-\chi_{0}(p)p^{-\alpha}}\right|^{-1} \nonumber \\
& \leq & \sup_{|t| \leq (yq)^{1/4}} e^{-\sum_{\sqrt{y} < p \leq y, p \nmid q} (1-\Re(\chi(p)p^{-it}))/p^{\alpha}} \nonumber \\
& \ll & e^{-y^{(1-\alpha)/2}/5} \nonumber \\
& \ll & e^{-\Theta(\sqrt{u\log u})}. \nonumber
\end{eqnarray}
We conclude that, under the conditions of Rodosski\u{\i} Bound 1 (and for ``large'' y),
$$ \left|\int_{\alpha - i(yq)^{1/4}}^{\alpha + i(yq)^{1/4}} L(s,\chi;y) x^{s} \breve{\Phi}(s) ds \right| \ll \sqrt{u\log u}e^{-\Theta(\sqrt{u\log u})} \Psi(x,y;\chi_{0},\Phi). $$
This estimate more than suffices for $\chi \in \Xi(k)$ with $k$ in the ``Rodosski\u{\i} range''.

\subsection{The remaining characters}
It remains to analyse $\Psi(x,y;\chi,\Phi)$ when $\chi \in \mathcal{A} := \bigcup_{k < 4A\log A + D} \Xi(k)$. Soundararajan's arguments work reasonably for this, the only adaptations being of the kind that we have demonstrated in $\S\S 2.3-2.4$, so a brief discussion seems sufficient.

We write $B=B(A):=\#\mathcal{A}$, which is bounded in terms of $A$ and $D$ because of the log-free zero density estimate (recall $\S 2.1$). The following result (which is proved by considering the corresponding sums with $\chi(p)p^{-it}$ replaced by $\chi(p)^{k}p^{-itk}$, for $1 \leq k \leq B+1$) is the content of Proposition 5.1 of Soundararajan~\cite{sound}:
\begin{rodbnd2}[Soundararajan, 2008]
Suppose that $\chi$ is a character of order exceeding $B=B(A)$. If $q \leq y^{A}$, and $|t| \leq q/(2(B+1))$, and $y/(A+1)^{2}$ is large, then
$$ \sum_{\sqrt{y} \leq p \leq y, p \nmid q} \frac{1 - \Re(\chi(p)p^{-it})}{p} \log p \geq \frac{\log y}{5(B+1)^{2}}. $$
\end{rodbnd2}
Using Rodosski\u{\i} Bound 2 in place of Rodosski\u{\i} Bound 1, one can proceed as in $\S 2.4$ to bound $\Psi(x,y;\chi,\Phi)$ for characters $\chi$ of order at least $B+1$ (with $e^{-\Theta(\sqrt{u\log u}/(B+1)^{2})}$ ultimately replacing $e^{-\Theta(\sqrt{u\log u})}$ in the estimates).

It now remains to treat $\chi \in \mathcal{B}$, where $\mathcal{B}:= \mathcal{A} \cap \{\chi: \chi \textrm{ has order } \leq B\}$. To simplify our formulae, we temporarily set $g=\log u /\sqrt{\log x \log y} = \log u /(\sqrt{u} \log y)$. A very small adaptation of the argument in $\S 2.4$ yields
\begin{eqnarray}
&& \left|\int_{\alpha - i(yq)^{1/4}}^{\alpha + i(yq)^{1/4}} L(s,\chi;y) x^{s} \breve{\Phi}(s) ds \right| \nonumber \\
& \ll & \left|\int_{\alpha - ig}^{\alpha + ig} L(s,\chi;y) x^{s} \breve{\Phi}(s) ds \right| + \sqrt{u\log u}\Psi(x,y;\chi_{0},\Phi) \sup_{g < |t| \leq (yq)^{1/4}} e^{-\sum_{\sqrt{y} < p \leq y} \frac{(1-\Re(\chi(p)p^{-it}))}{p^{\alpha}}} \nonumber \\
& \ll & \frac{x^{\alpha} \breve{\Phi}(\alpha) \log u}{\sqrt{\log x \log y}} L(\alpha,\chi_{0};y) \sup_{|t| \leq g} e^{-\sum_{p \leq y, p \nmid q} \frac{(1-\Re(\chi(p)p^{-it}))}{p^{\alpha}}} + \nonumber \\
&& + \sqrt{u\log u}\Psi(x,y;\chi_{0},\Phi) \sup_{g < |t| \leq (yq)^{1/4}} e^{-\sum_{\sqrt{y} < p \leq y, p \nmid q} \frac{(1-\Re(\chi(p)p^{-it}))}{p^{\alpha}}}, \nonumber
\end{eqnarray}
if $y$ is ``large''. By the argument of Lemma 5.2 of Soundararajan~\cite{sound}, (which is proved by a neat reduction from working with characters of order at most $B(A)$ to working with the principal character), provided that $u \geq B^{3}$ (say) the exponential in the second term is $ \ll e^{-\Theta(\log^{2}u)}$. This gives an acceptable bound for that term.

Finally we apply the following result, which may be extracted\footnote{It requires a little care to obtain the explicit error terms in Character Sum Bound 1. The reader should note that if $\chi$ is a (primitive) non-principal character to modulus $q$, of order $l$, and $\delta > 0$ is small, then Heath-Brown's~\cite{hb} refined character sum estimate shows in particular that
$$ \sum_{n \leq H} \chi(n) \ll_{\delta} H l^{3\delta/2} (q^{1/4}H^{-1})^{\delta} q^{\delta^{2}/2}. $$
This should be applied with $H = z/d \geq y^{\sqrt{e}-3\delta/2} \geq q^{1/4}y^{7\delta/2}$ in the proof of Soundararajan's~\cite{sound} Lemma 5.3.} from Lemmas 5.2 and 5.3 of Soundararajan~\cite{sound}.
\begin{charbnd}[Soundararajan, 2008]
Let $\delta > 0$ be sufficiently small, and suppose that $q \leq y^{4\sqrt{e}-20\delta}$. If $\chi \neq \chi_{0}$ is a character of order at most $B$, and if $|t| \leq 1/(B\log y)$, then
$$ \sum_{p \leq y, p \nmid q} \frac{1-\Re(\chi(p)p^{-it})}{p^{\alpha}} \gg \frac{1}{B} \log((\delta+O(\frac{1}{\log y})+O_{\delta}(\frac{B^{3\delta/2} \log^{2}y}{y^{\delta^{2}}})) y^{\delta(1-\alpha)/2}) + O(1). $$
\end{charbnd}
Using this to bound the exponential in the first term, we complete the proof of Theorem 1 for ``large'' values of $y$. Character Sum Bound 1 imports Burgess's character sum estimates (with the modification of Heath-Brown~\cite{hb} for characters of bounded order): indeed, it is clear that when $|t|$ is so small, cancellation in $L(\alpha+it,\chi;y)$ amounts to cancellation in sums of $\chi$. If one could prove non-trivial bounds for shorter character sums, one could introduce them at this point and thereby extend the range of $q$ in Theorem 1. These remarks apply equally to Soundararajan's proof~\cite{sound} on the range of $y$ where it is valid.

\subsection{The case of small $y$}
When $y$ is ``small'' the argument from Soundararajan's paper~\cite{sound} almost goes through, even in a simplified form. The complication is that the parameter $\alpha(x,y)$ behaves differently when $y \leq \log x$ than otherwise, and in particular is very small, which necessitates a few changes. We sketch these briefly.

Thus when $y \leq \log x$, we have
\begin{eqnarray}
\left|\frac{L(\alpha+it,\chi;y)}{L(\alpha,\chi_{0};y)}\right| = \prod_{p \leq y, p \nmid q} \left|1+\frac{1-\chi(p)p^{-it}}{p^{\alpha}-1}\right|^{-1} & \leq & \prod_{p \leq y, p \nmid q} \left|1+\frac{1-\Re(\chi(p)p^{-it})}{y^{\alpha}-1}\right|^{-1} \nonumber \\
& \leq & \prod_{p \leq y, p \nmid q} \left|1+\frac{c\log x (1-\Re(\chi(p)p^{-it}))}{y}\right|^{-1} \nonumber
\end{eqnarray}
for a small constant $c > 0$, since $\alpha(x,y) = O(y/\log x \log y)$ for $y \leq \log x$. By calculus, $1+ct \geq (1+t)^{c}$ for $t \geq 0$ and $0 \leq c \leq 1$, so the above is at most
$$ \left(1+\frac{\log x}{y} \right)^{-c\sum_{p \leq y, p \nmid q} (1-\Re(\chi(p)p^{-it}))} \leq \left(1+\frac{\log x}{y} \right)^{-c \frac{\sqrt{y}}{\log y} \sum_{\sqrt{y} \leq p \leq y, p \nmid q} \frac{1-\Re(\chi(p)p^{-it})}{p} \log p}. $$
It follows, using the decay of $\breve{\Phi}$ to control the large range of integration, that
\begin{eqnarray}
&& \left|\int_{\alpha - i(yq)^{1/4}}^{\alpha + i(yq)^{1/4}} L(s,\chi;y) x^{s} \breve{\Phi}(s) ds \right| \nonumber \\
& \ll & x^{\alpha} \breve{\Phi}(\alpha) L(\alpha,\chi_{0};y) \sup_{|t| \leq (yq)^{1/4}} \left(1+\frac{\log x}{y} \right)^{-c \frac{\sqrt{y}}{\log y} \sum_{\sqrt{y} \leq p \leq y, p \nmid q} \frac{1-\Re(\chi(p)p^{-it})}{p} \log p}. \nonumber
\end{eqnarray}
If $\chi \in \Xi(k)$ for some $k \geq 4A\log A + D$ then Rodosski\u{\i} Bound 1 shows the above is $O(\Psi(x,y;\chi_{0},\Phi)/(\log x + 2^{y^{1/3}}))$, say (since we assume $1 \ll y \leq \log x$).

When $\log x < y \leq e^{\log^{1/10}x}$, one can argue instead (as in $\S 2.4$) that
$$ \left|\frac{L(\alpha+it,\chi;y)}{L(\alpha,\chi_{0};y)}\right| \leq \prod_{p \leq y, p \nmid q} \left(1+\sum_{k=1}^{\infty} \frac{1-\Re(\chi(p)p^{-it})}{p^{k\alpha}}\right)^{-1} \leq e^{-\frac{y^{(1-\alpha)/2}}{\log y} \sum_{\sqrt{y} \leq p \leq y, p \nmid q} \frac{1-\Re(\chi(p)p^{-it})}{p} \log p}, $$
so Rodosski\u{\i} Bound 1 implies that the integral is $O(\Psi(x,y;\chi_{0},\Phi)/2^{\log^{4}y})$, say (as $y^{1-\alpha} \gg u \geq \log^{9/10}x \geq \log^{9}y$). These estimates are more than acceptable for $\chi \in \Xi(k)$, $4A\log A + D \leq k \leq (\log q)/2$.

We must still deal with $\chi \in \Xi(k)$ for $k$ in the ``problem range''. As in $\S 2.5$, if $\chi$ has order at least $B(A)+1$ then we need only apply Rodosski\u{\i} Bound 2 instead of Rodosski\u{\i} Bound 1 in the preceding calculations. If $\chi \in \mathcal{B}$ we set $h=\log u /\sqrt{\min\{\log x \log y , y \log y\}}$, so $|\int_{\alpha - i(yq)^{1/4}}^{\alpha + i(yq)^{1/4}} L(s,\chi;y) x^{s} \breve{\Phi}(s) ds|$ is
$$ \ll x^{\alpha} \breve{\Phi}(\alpha) (h \sup_{|t| \leq h} |L(\alpha+it,\chi;y)| + L(\alpha,\chi_{0};y) \sup_{h < |t| \leq (yq)^{1/4}} e^{-c\sum_{p \leq y, p \nmid q} \frac{(1-\Re(\chi(p)p^{-it}))}{p^{\alpha}}}). $$
Provided that $u \geq B^{3}$ and $y \geq B^{7}$ (say), the argument of Lemma 5.2 of Soundararajan~\cite{sound} shows the exponential is $\ll e^{-\Theta(\log^{2}u)} + e^{-\Theta(y^{2/3})}$, which is  $\ll e^{-\Theta((\log\log x)^{2})}$ for ``small'' $y$ (that is $(\log\log x)^{3} \leq y \leq e^{\log^{1/10}x}$). This is also acceptable.

Finally we apply Character Sum Bound 1 to estimate $\sup_{|t| \leq h} |L(\alpha+it,\chi;y)|$. If $y \geq \log x$ this works precisely as in $\S 2.5$. If $(\log\log x)^{3} \leq y \leq \log x$ we have
$$ \sup_{|t| \leq h} |L(\alpha+it,\chi;y)| \leq \left(1+\frac{\log x}{y} \right)^{-(c/B) \log((\delta+o(1)) y^{\delta(1-\alpha)/2}) + O(1)} L(\alpha,\chi_{0};y), $$
noting that $h \leq 1/\log^{2}y$ (say) for such $y$, so Character Sum Bound 1 is applicable. Now $x^{\alpha} \breve{\Phi}(\alpha) h L(\alpha,\chi_{0};y) \ll (\Psi(x,y;\chi_{0},\Phi) \log u \log x) / y$ when $y \leq \log x$, and that is $\ll \Psi(x,y;\chi_{0},\Phi) \log y ((\log x)/y)^{2}$. This all implies an acceptable bound for $x^{\alpha} \breve{\Phi}(\alpha) h \sup_{|t| \leq h} |L(\alpha+it,\chi;y)|$, finishing the proof of Theorem 1 for ``small'' $y$.

\subsection{The case of very small $y$}
If $y$ is ``very small'', two changes are needed to the argument for ``small'' $y$ in $\S 2.6$ (and at the end of $\S 2.1$).

Firstly, if $q\sqrt{y} \leq (\log\log x)^{1/3}$ we cannot bound integrals like $\int_{\alpha + i(yq)^{1/4}}^{\alpha + i\infty} L(s,\chi;y) x^{s} \breve{\Phi}(s) ds$ acceptably just using the decay of $\breve{\Phi}$, necessarily. As an alternative, if $y \leq \sqrt{\log x}$, if $p_{0} \leq y$ is any prime not dividing $q$, and if $\epsilon > 0$, then as in $\S 2.6$ we have
\begin{eqnarray}
&& \left|\int_{|t| \geq \epsilon} L(\alpha+it,\chi;y) x^{\alpha+it} \breve{\Phi}(\alpha+it) dt \right| \nonumber \\
& \ll & x^{\alpha} L(\alpha,\chi_{0};y) \int_{|t| \geq \epsilon} \left(1+\frac{\log x}{y} \right)^{-c(1-\Re(\chi(p_{0})p_{0}^{-it}))} \frac{1}{|t|(|t|+1)^{8}} dt \nonumber \\
& \ll & \Psi(x,y;\chi_{0},\Phi) \sqrt{\frac{y}{\log y}} \int_{|t| \geq \epsilon} (\log x)^{-c(1-\Re(\chi(p_{0})p_{0}^{-it}))/2} \frac{1}{|t|(|t|+1)^{8}} dt, \nonumber
\end{eqnarray}
on recalling the lower bound for $\Psi(x,y;\chi_{0},\Phi)$ (and $\breve{\Phi}(\alpha)$) from $\S 2.1$. But
$$ 1-\Re(\chi(p_{0})p_{0}^{-it}) = 1 - \cos(\arg(\chi(p_{0}))-t\log p_{0}), $$
which is clearly $\geq 1/(\log\log x)^{4/5}$, say, except on a progression of intervals of $t$ having lengths $\Theta(1/(\log p_{0} (\log\log x)^{2/5}))$ and spacing $\Theta(1/\log p_{0})$. Thus if $\epsilon \geq 1/(\log\log x)^{2/5-1/3}$, the integral is
$$ \ll \Psi(x,y;\chi_{0},\Phi) \sqrt{\frac{y}{\log y}} \frac{1}{\epsilon} (\frac{1}{(\log\log x)^{2/5}} + e^{-\Theta((\log\log x)^{1/5})}) \ll \frac{\sqrt{y} \Psi(x,y;\chi_{0},\Phi)}{(\log\log x)^{1/3}}. $$
Note that this holds for {\em all} Dirichlet characters $\chi$ to modulus $q$. If $q\sqrt{y} \leq (\log\log x)^{1/3}$ then we apply this with $\epsilon$ chosen as $\min\{q/(2(B(A)+1)),1\}$, (which is at least $1/y^{2/5-1/3}$ by assumption that $y$ is large in terms of $A$, and thus at least $1/(\log\log x)^{2/5-1/3}$). Then we can estimate $\int_{|t| \leq \epsilon} L(\alpha+it,\chi;y) x^{\alpha+it} \breve{\Phi}(\alpha+it) dt $ for $\chi \notin \mathcal{B}$ using the Rodosski\u{\i} Bounds, exactly as demonstrated in $\S 2.6$.

(Note that we need not assume that $q \geq \sqrt{y}$, as previously, for the Rodosski\u{\i} Bounds to apply, since we are concerned with $|t| \leq \epsilon \leq q/(2(B+1))$ rather than $|t| \leq (yq)^{1/4}$.)

Secondly, to deal with $\chi \in \mathcal{B}$ we just set $h=1/\log^{2}y$, rather than choosing $h$ as in $\S 2.6$. Then $|\int_{\alpha - i(yq)^{1/4}}^{\alpha + i(yq)^{1/4}} L(s,\chi;y) x^{s} \breve{\Phi}(s) ds|$, or $|\int_{\alpha - i\epsilon}^{\alpha + i\epsilon} L(s,\chi;y) x^{s} \breve{\Phi}(s) ds|$, is
$$ \ll x^{\alpha} \breve{\Phi}(\alpha) (h \sup_{|t| \leq h} |L(\alpha+it,\chi;y)| + L(\alpha,\chi_{0};y) \sup_{h < |t| \leq (yq)^{1/4}} (1+\frac{\log x}{y})^{-c\sum_{p \leq y, p \nmid q} \frac{(1-\Re(\chi(p)p^{-it}))}{p^{\alpha}}}), $$
where these terms may be bounded as in $\S 2.6$. In particular the exponential is $\ll (\log x)^{-\Theta(y^{2/3})}$ if $B^{7} \leq y \leq (\log\log x)^{3}$, which is more than satisfactory.
\begin{flushright}
Q.E.D.
\end{flushright}

\section{Proof of Proposition 1}
To start the proof of Proposition 1, as used in $\S 2.3$, we shall establish the following lemma. We will need the result, and the techniques of the proof will also be used again.
\begin{lem3}
Under the assumptions of Proposition 1, and for primitive $\chi$, we have
$$ \left|\frac{L'(\sigma+it,\chi)}{L(\sigma+it,\chi)} \right| \ll \log q. $$
\end{lem3}
To see this, we note that the left hand side is certainly at most
$$ \left|\frac{L'(1+1/\log q+it,\chi)}{L(1+1/\log q+it,\chi)}\right| + (1+\frac{1}{\log q}-\sigma) \sup_{\sigma \leq \sigma' \leq 1+1/\log q} \left|\frac{d}{d \sigma'} \frac{L'(\sigma'+it,\chi)}{L(\sigma'+it,\chi)} \right|. $$
Here the first term is at most $\zeta'(1+1/\log q)/\zeta(1+1/\log q)$, which is $O(\log q)$. We also note that
$$1 + 1/\log q - \sigma \leq (1-\alpha(x,y)) + 1/\log q + Bk/\log x \leq k/(4\log q), $$
say, in view of the assumptions on $y,q$ and $k$ in Proposition 1.

For primitive $\chi$, differentiation of the explicit formula for $\frac{L'(s,\chi)}{L(s,\chi)}$ (which is e.g. formula (17) in chapter 12 of Davenport~\cite{davenport}) yields
$$ \frac{d}{d\sigma'} \frac{L'(\sigma'+it,\chi)}{L(\sigma'+it,\chi)} = -\sum_{n=0}^{\infty} \frac{1}{(2n+\sigma'+it+a(\chi))^{2}} - \sum_{\rho} \frac{1}{(\sigma'+it-\rho)^{2}}, $$
where the second sum is over the non-trivial zeros of $L(s,\chi)$, and $a(\chi)$ is 0 or 1 according as $\chi(-1)$ is 1 or -1. (Thus the first sum is really over the trivial zeros of $L(s,\chi)$: see e.g. chapters 9 and 19 of Davenport~\cite{davenport}). The sum over $n$ is clearly $O(1)$, and since $\sigma' \geq \sigma \geq 1-k/(4\log q)$, $|t| \leq q/2$ and $\chi \in \Xi(k)$ we have
\begin{eqnarray}
\sum_{\rho} \frac{1}{|\sigma'+it-\rho|^{2}} & \leq & \sum_{|\Im(\rho)-t| \leq 1, \atop |\Im(\rho)| \leq q} \frac{1}{|\sigma'+it-\rho|^{2}} + \sum_{|\Im(\rho)-t| > 1, \atop |\Im(\rho)| \leq q} \frac{1}{|\sigma'+it-\rho|^{2}} + \sum_{|\Im(\rho)| > q} \frac{4}{|\Im(\rho)|^{2}} \nonumber \\
& \ll & \sum_{|\Im(\rho)-t| \leq 1} \frac{1}{|1+1/\log q+it-\rho|^{2}} + \log q + \frac{\log q}{q}, \nonumber
\end{eqnarray}
using the fact that $\Re(\rho) \leq 1 - k/\log q$ in the first sum, and standard results on the vertical distribution of zeros of $L(s,\chi)$ (as in e.g. chapter 16 of Davenport~\cite{davenport}).

To bound the remaining sum, we again use the fact that $\chi \in \Xi(k)$, noting that
\begin{eqnarray}
\sum_{|\Im(\rho)-t| \leq 1} \frac{1}{|1+1/\log q + it-\rho|^{2}} & = & \sum_{|\Im(\rho)-t| \leq 1} \frac{1}{\Re(1+1/\log q + it - \rho)} \Re(\frac{1}{1+1/\log q + it -\rho}) \nonumber \\
& \leq & \frac{\log q}{k} \Re(\sum_{|\Im(\rho)-t| \leq 1} \frac{1}{1+1/\log q + it -\rho}) \nonumber \\
& = & \frac{\log q}{k} \Re(\frac{L'(1+1/\log q + it,\chi)}{L(1+1/\log q + it,\chi)} + O(\log q)) \nonumber \\
& \ll & \frac{\log^{2}q}{k}. \nonumber
\end{eqnarray}
Here the second equality is a classical approximation for $L'(s,\chi)/L(s,\chi)$, as in formula (4) of chapter 16 of Davenport~\cite{davenport}. See section 4 of Soundararajan~\cite{sound}, and especially chapter 9.2 of Montgomery~\cite{mont}, for further illustration of this argument.

Combining the estimates we obtained proves Lemma 3.
\vspace{12pt}

Now we note that, under the conditions of Proposition 1,
\begin{eqnarray}
|\log L(\sigma+it,\chi;y) - \log L(\alpha+it,\chi;y)| & \leq & (\alpha-\sigma) \sup_{\sigma \leq \sigma' \leq \alpha} \left|\frac{L'(\sigma'+it,\chi;y)}{L(\sigma'+it,\chi;y)} \right| \nonumber \\
& \leq & \frac{Bk}{\log x} (\sup_{\alpha-\frac{Bk}{\log x} \leq \sigma' \leq \alpha} \left|\sum_{n \leq y} \frac{\Lambda(n) \chi(n)}{n^{\sigma'+it}} \right| + O(1)). \nonumber
\end{eqnarray}
We suppose initially that $\chi$ is a primitive Dirichlet character, and also that it will suffice to bound the above with the sum over $n$ replaced by $\sum_{n \leq Ry} w(n) \frac{\Lambda(n)\chi(n)}{n^{\sigma'+it}}$, where $R:=\max\{2,y^{y^{-k/2\log q}}\}$ and
$$ w(n) := \left\{ \begin{array}{ll}
     1 & \textrm{if } 1 \leq n \leq y   \\
     1 - \frac{\log(n/y)}{\log R} & \textrm{if } y \leq n \leq Ry   \\
     0 & \textrm{otherwise}
\end{array} \right .$$
At the end of the argument we will show how to remove these assumptions.

Recall that for $\Re(s)>1$ we have $\sum_{n=1}^{\infty} \frac{\Lambda(n) \chi(n)}{n^{s}} = - \frac{L'(s,\chi)}{L(s,\chi)}$. Then a fairly standard contour integration procedure, as in e.g. chapters 13.2 and 12.1 of Montgomery and Vaughan~\cite{mv}, reveals that
\begin{eqnarray}
\sum_{n \leq Ry} w(n) \frac{\Lambda(n) \chi(n)}{n^{s}} & = & -\frac{L'(s,\chi)}{L(s,\chi)} - \frac{1}{\log R} \sum_{n=0}^{\infty} \frac{(Ry)^{-2n-a(\chi)-s}-y^{-2n-a(\chi)-s}}{(2n+a(\chi)+s)^{2}} - \nonumber \\
&& - \frac{1}{\log R} \sum_{\rho} \frac{(Ry)^{\rho-s} - y^{\rho-s}}{(\rho-s)^{2}} \nonumber
\end{eqnarray}
whenever $L(s,\chi) \neq 0$. Here our notation is exactly as above. As always, the purpose of introducing the smoother weight $w(n)$ was to obtain nicer behaviour in these sums, namely that all denominators are raised to at least the second power. This is also the reason that it was a good idea to switch to studying $\frac{d}{d \sigma'} \frac{L'(\sigma'+it,\chi)}{L(\sigma'+it,\chi)}$ in proving Lemma 3.

Putting $s=\sigma'+it$, Lemma 3 and a trivial estimation show the first two terms in the above are $O(\log q)$. To estimate the sum over $\rho$ we proceed as in the proof of Lemma 3, noting that we can extract a power saving $O(y^{-3k/(4\log q)})$ on the range $|\Im(\rho)| \leq q$ (since $\sigma' \geq 1-k/(4\log q)$ and $\chi \in \Xi(k)$). Thus we have
\begin{eqnarray}
\left|\frac{1}{\log R} \sum_{\rho} \frac{(Ry)^{\rho-\sigma'-it} - y^{\rho-\sigma'-it}}{(\rho-\sigma'-it)^{2}} \right| & \ll & \frac{y^{-3k/(4\log q)} \log^{2}q}{k \log R} + \frac{1}{\log R} \sum_{|\Im(\rho)| > q} \frac{(Ry)^{k/(4\log q)}}{|\Im(\rho)|^{2}} \nonumber \\
& \ll & \frac{A\log q}{k} + \frac{y^{1/4}\log q}{q}, \nonumber
\end{eqnarray}
since $\log R \geq y^{-k/(2\log q)} \log y$ but $Ry \leq y^{2}$. This is acceptable for Proposition 1.

It remains to justify the two assumptions that we made at the start of the proof. Firstly, if $k$ is such that $R=y^{y^{-k/2\log q}} \geq 2$ then
$$ \left|\sum_{y < n \leq Ry} w(n) \frac{\Lambda(n)\chi(n)}{n^{\sigma'+it}} \right| \leq \sum_{y < n \leq Ry} \frac{\Lambda(n)}{n^{\sigma'}} \ll \frac{y^{1-\sigma'}(R^{1-\sigma'}-1)}{1-\sigma'}, $$
and since $1 - \sigma' \leq k/(4\log q)$ this is $\ll y^{1-\sigma'}\log R$, which is at most $\log y$. If $k$ is such that $R=2$, the same argument produces a bound $O(y^{1-\sigma'})$, or alternatively
$$ \left|\sum_{y < n \leq 2y} w(n) \frac{\Lambda(n)\chi(n)}{n^{\sigma'+it}} \right| \leq \frac{1}{y^{\sigma'}} \max_{y < m \leq 2y} \left|\sum_{y < n \leq m} \frac{\Lambda(n)\chi(n)}{n^{it}} \right| \ll y^{1-\sigma'-k/\log q} \log^{2}q. $$
Here the first inequality follows from Abel's partial summation lemma, and the second from Lemma 3.1 of Soundararajan~\cite{sound} (or directly by an explicit formula argument rather easier than the above calculations). Comparing our two bounds, we see that when $R=2$ the sum must be $\ll y^{1-\sigma'-k/2\log q}\log q$, which is at most $\log q$. These error estimates are acceptable for Proposition 1.

Finally, if $\chi$ is not primitive we can apply the above techniques to the primitive character inducing $\chi$. This results in an error term of size at most
$$ \sum_{p|q} \log p \sum_{r \geq 1} \frac{1}{p^{r\sigma'}} \ll \sum_{p|q} \log p \leq \log q $$
when estimating $\sum_{n \leq y} \frac{\Lambda(n) \chi(n)}{n^{\sigma'+it}}$, which again is acceptable for Proposition 1.
\begin{flushright}
Q.E.D.
\end{flushright}

\appendix
\section{Unsmoothing}
In this appendix we briefly explain how to pass from results about $\Psi(x,y;q,a,\Phi)$, which we actually proved, to results about the unsmoothed version $\Psi(x,y;q,a)$.

To ``unsmooth'' one notes that if $\epsilon > 0$, and $\textbf{1}_{[0,1-\epsilon]} \leq \Phi \leq \textbf{1}_{[0,1]}$, then
\begin{eqnarray}
\Psi(x,y;q,a) \geq \Psi(x,y;q,a,\Phi) & = &  \frac{1}{\phi(q)}\Psi(x,y;\chi_{0},\Phi)(1+o_{\Phi}(1)) \nonumber \\
& \geq & \frac{1}{\phi(q)}\Psi_{q}(x,y)(1-\frac{(\Psi_{q}(x,y)-\Psi_{q}((1-\epsilon)x,y))}{\Psi_{q}(x,y)} + o_{\Phi}(1)), \nonumber
\end{eqnarray}
where the first equality is what we proved in the body of this paper. One can similarly obtain an upper bound for $\Psi(x,y;q,a)$, so to deduce Theorem 1 we need to know that for any $\eta > 0$, the ratio $(\Psi_{q}(x,y)-\Psi_{q}((1-\epsilon)x,y))/\Psi_{q}(x,y)$ will be at most $\eta$ if $\epsilon$ is chosen sufficiently small (and $\log x/\log q$ is large enough). This local result about $\Psi_{q}(x,y)$ follows from Th\'{e}or\`{e}me 2.4 of de la Bret\`{e}che and Tenenbaum~\cite{dlbten}, (also see Theorem 3 of Hildebrand and Tenenbaum~\cite{ht}), except when $y$ does not tend to infinity with $x$. However, if $2 \leq y \leq \sqrt{\log x}$ one has
$$ \Psi_{q}(x,y) = \frac{1}{(\#\{p \leq y : p \textrm{ is prime}, p \nmid q\})!} \prod_{p \leq y, p \nmid q} \left(\frac{\log x}{\log p}\right) \left(1+O\left(\frac{y^{2}}{\log x \log y}\right) \right), $$  
exactly similarly to an expression for $\Psi(x,y)$ due to Ennola (and explained in Chapter III.5.2 of Tenenbaum's book~\cite{ten}), which directly implies that $\Psi_{q}((1-\epsilon)x,y)$ is $(1+O(\epsilon))\Psi_{q}(x,y)$, say. Actually one obtains a bound $O(\epsilon)$ for our ratio in all cases, provided that $\log x / \log q$ is large enough in terms of $\epsilon$.

For the proof of Theorem 2, one should apply this procedure to $\Psi(x,y;\chi,\Phi)$ for all of the (bounded number of) characters $\chi \in \mathcal{B}$. The ``analytic'' unsmoothing procedure used by Soundararajan~\cite{sound} appears not to work on our extended range of $y$, because one cannot replace $\Psi(x,y;\chi,\Phi)$ by a contour integral over a suitably short range of $t$ for $\breve{\Phi}$ to be removed from it.

\section{Estimates for Dirichlet series involving the principal character}
In this appendix we prove two estimates for Dirichlet series involving the principal character, and a bound for an oscillating integral, which were needed in $\S\S 2.3-2.5$. These results are of a rather standard type (see e.g. Lemma 8 of Hildebrand and Tenenbaum~\cite{ht}), but we include the short proofs in the interests of completeness.

We suppose that $3/4 \leq \beta \leq 1$, say: in $\S\S 2.3-2.5$ we had $\beta = \alpha(x,y)$ or $\beta = \alpha(x,y) - Bk/\log x$. We also suppose that $\log q \leq \log^{2}y$, which certainly implies that $ \sum_{p \geq \log y, p \mid q} 1/p \ll 1$. By partial summation from the prime number theorem, if $2 \leq z \leq y$ and if $t \neq 0$ then
$$ \sum_{z \leq p \leq y} \frac{1-\cos(t\log p)}{p} = \log\log y - \log\log z - \int_{t\log z}^{t\log y} \frac{\cos w}{w} dw + O((1+|t|)e^{-d\sqrt{\log z}}), $$
for a certain constant $d > 0$. Choosing $z=e^{(\log\log y)^{3}}$, we find that if $1/\log y \leq |t| \leq 1/(\log\log y)^{3}$ then the sum is $\log\log y + \log |t| + O(1)$, whilst if $1/(\log\log y)^{3} \leq |t| \leq \log^{1/4}y$, say, then the sum is $\log\log y - 3\log\log\log y + O(1)$. Thus we have
\begin{eqnarray}
\left|\frac{L(\beta+it,\chi_{0};y)}{L(\beta,\chi_{0};y)}\right| \leq e^{-\sum_{p \leq y, p \nmid q} (1-\cos(t\log p))/p} & \ll & e^{-\sum_{\log y \leq p \leq y} (1-\cos(t\log p))/p} \nonumber \\
& \ll & \frac{\max\{|t|^{-1},(\log\log y)^{3}\}}{\log y}, \nonumber
\end{eqnarray}
provided that $1/\log y \leq |t| \leq \log^{1/4}y$, and so
$$ \int_{\beta-i\log^{1/4}y}^{\beta+i\log^{1/4}y} \left|L(s,\chi_{0};y)\right|^{2} d|s| = O(|L(\beta,\chi_{0};y)|^{2}/\log y). $$

Next, we note that
$$ \left|\frac{L'(\beta+it,1;y)}{L(\beta+it,1;y)}\right| = \left|\sum_{p \leq y} \frac{\log p}{p^{\beta+it}-1} \right| = \left|\sum_{n \leq y} \frac{\Lambda(n)}{n^{\beta+it}}\right| + O(1), $$
since $\beta$ is large. By partial summation from the prime number theorem, if $y \geq 2$ then
$$ \sum_{n \leq y} \frac{\Lambda(n)}{n^{\beta+it}} = \frac{y^{1-\beta-it}-1}{1-\beta-it} + O((1+|t|)\int_{1}^{y} \frac{e^{-d\sqrt{\log w}}}{w^{\beta}} dw) = \frac{y^{1-\beta-it}-1}{1-\beta-it} + O((1+|t|)y^{1-\beta}), $$
so if $\beta \leq \alpha = \alpha(x,y)$ then certainly
$$ \sum_{n \leq y} \frac{\Lambda(n)}{n^{\beta+it}} \ll y^{\alpha-\beta} u \min\{\frac{\log u}{|t|}, \log y\} + (1+|t|)y^{\alpha-\beta} u \log u, $$
on recalling the definition of $\alpha$. It follows, provided $\log u \leq \log^{1/4}y$ (which certainly holds if $y$ is ``large'', in the sense of the introduction), that
$$ \int_{\beta-i\log^{1/4}y}^{\beta+i\log^{1/4}y} \left|\frac{L'(s,1;y)}{L(s,1;y)}\right|^{2} d|s| = O(y^{2(\alpha-\beta)} u^{2}\log u \log y). $$

Finally note that if measurable $0 \leq \Phi \leq 1$ is supported on $[0,2]$, then
$$ \int_{t}^{\log^{1/4}y} x^{ir} \breve{\Phi}(\beta+ir) dr = \int_{0}^{2} \Phi(v) v^{\beta-1} \int_{t}^{\log^{1/4}y} (xv)^{ir} dr dv , $$
by definition of the Mellin transform. Splitting the integral over $v$ at $1/\sqrt{x}$, and evaluating the integral over $r$, we find that this is $O(1/\log x)$.

\vspace{12pt}
\noindent {\em Acknowledgements.} The author would like to thank Andrew Granville, Ben Green and Kannan Soundararajan for their encouragement, and for comments on a draft of this paper.



\begin{thebibliography}{99}

\bibitem{davenport} H. Davenport. {\em Multiplicative Number Theory.} Third edition, revised by H. L. Montgomery, published by Springer. 2000.

\bibitem{dlbten} R. de la Bret\`{e}che, G. Tenenbaum. Propri\'{e}t\'{e}s statistiques des entiers friables. {\em The Ramanujan Journal}, \textbf{9}, pp 139-202. 2005

\bibitem{granville} A. Granville. Integers, without large prime factors, in arithmetic progressions. II {\em Phil. Trans. R. Soc. Lond. A}, \textbf{345}, pp 349-362. 1993

\bibitem{hb} D. R. Heath-Brown. Zero-free regions for Dirichlet $L$-functions, and the least prime in an arithmetic progression. {\em Proc. London Math. Soc.}, \textbf{64}, pp 265-338. 1992

\bibitem{ht} A. Hildebrand, G. Tenenbaum. On integers free of large prime factors. {\em Trans. Amer. Math. Soc.,} \textbf{296}, no. 1, pp 265-290. 1986

\bibitem{mont} H. L. Montgomery. {\em Ten Lectures on the Interface Between Analytic Number Theory and Harmonic Analysis.} Published for the Conference Board of the Mathematical Sciences by the American Mathematical Society. 1994

\bibitem{mv} H. L. Montgomery, R. C. Vaughan. {\em Multiplicative Number Theory I: Classical Theory.} First edition, published by Cambridge University Press. 2007

\bibitem{sound} K. Soundararajan. The distribution of smooth numbers in arithmetic progressions. {\em Anatomy of Integers,} CRM Proc. and Lect. Notes, vol. 46, Amer. Math. Soc., Providence, RI, pp 115-128. 2008

\bibitem{ten} G. Tenenbaum. {\em Introduction to analytic and probabilistic number theory.} English edition, published by Cambridge University Press. 1995

\end{thebibliography}
\end{document}